# Différences entre approches semi-quantitative et quantitative pour l'évaluation probabiliste des risques technologiques

# Differences between semi-quantitative and quantitative approaches for probabilistic assessment of technological risks

Elsa Rosner, Florent Brissaud, Bruno Declerck
DNV France
69 rue du Chevaleret
75013 Paris
+33 (0)1 44 24 40 13
elsa.rosner@dnvgl.com

Yann Flauw, Valérie de Dianous
INERIS
BP 2
60550 Verneuil-en-Halatte

**Résumé**

L'objectif de cette communication est de comparer l'approche quantitative de type arbre de défaillances / arbre d'évènements à l'approche dite "semi-quantitative" telle que présentée par le guide Oméga 10 (INERIS, 2008), pour les démarches d'évaluation probabiliste des risques technologiques. Après une brève présentation de l'approche quantitative "classique" basée sur les arbres de défaillances et d'évènements, ainsi que de l'approche "semi-quantitative" du guide Oméga 10, un cas d'étude est présenté afin d'illustrer et de comparer leurs mises en application. Ce cas d'étude porte sur un séparateur (eau/gaz/condensats) basse pression pour lequel des scénarios ont été identifiés autour de l'évènement redouté "dégagement de gaz suite à une surpression dans le séparateur basse-pression". Outre une boucle de contrôle, ce système comprend les barrières de sécurité suivantes : une alarme avec action d'un opérateur, un système instrumenté de sécurité, deux soupapes, et un cloisonnement anti-explosion. Ce cas d'étude représente ainsi un système relativement conventionnel dans l'industrie. Pour ce qui est de la mise en œuvre de l'une et l'autre des approches, l'évaluation "semi-quantitative" selon le guide Oméga 10 est plus facile et plus rapide. Pour ce qui est des résultats, certains écarts sont observables, allant jusqu'à un facteur 10.

**Summary**

The aim of this paper is to compare the probabilistic approach based on fault trees and event trees, to an approach called "semi-quantitative" as presented by the guide Omega 10 (INERIS, 2008), for the industrial risk management. A brief presentation of the "classical" probabilistic risk assessment as well as the "semi-quantitative" approach is proposed. Then, a case study is detailed to illustrate and compare their use. This case study focuses on a low pressure separator (water/gas/condensate) for which scenarios have been identified around the undesired event "gas leak due to overpressure in the separator." The system includes a control loop and the following safety barriers: an alarm with operator actions, a safety instrumented system, two safety relief valves, and a blast-proof partition. This case study therefore represents a relatively conventional system in the industry. Regarding to the implementation of the two approaches, the assessment according to the Omega 10 guide is easier and faster. In terms of results, differences can be observed, up to a factor 10.

## 1. Introduction

En France, la réglementation sur les risques technologiques repose sur celle des Installations Classées pour la Protection de l'Environnement (ICPE), instaurée dans les années soixante-dix et aujourd'hui régie par le Livre V (Titre I) du Code de l'Environnement (Journal Officiel, 2000). À l'origine, ce cadre est "déterministe", c'est-à-dire qu'il ne prend en compte que certains scénarios d'accidents de référence (i.e. les "pires cas" (Lenoble et Durand, 2011)). Parallèlement à l'instauration des ICPE, les années soixante-dix voient aussi apparaître des approches "quantitatives" appliquées à l'estimation de risques dans le nucléaire aux USA, mais également en France (Lannoy, 2008). Ce type d'approche trouve aussi écho dans l'industrie pétrolière au Royaume-Uni, à la fin de cette même décennie (Lannoy, 2008). Pour l'industrie française en général, c'est à la suite de la catastrophe d'AZF à Toulouse, le 21 septembre 2001, que la réglementation sur les risques technologiques a intégré des critères probabilistes (Lenoble et Durand, 2011) avec la "loi Risques" (Chapitre II) du 30 juillet 2003 (Journal Officiel, 2003).

En particulier, la "loi Risques" a fait évoluer la réglementation sur les études de dangers (EDD). Les EDD sont requises lors des demandes d'autorisation des ICPE (ce qui concerne plus de 43 000 établissements en France). Le Code de l'Environnement, ainsi modifié par la "loi Risques" (Journal Officiel, 2003), mentionne que "cette étude donne lieu à une analyse de risques qui prend en compte la probabilité d'occurrence, la cinétique et la gravité des accidents potentiels selon une méthodologie qu'elle explicite" et qu' "elle définit et justifie les mesures propres à réduire la probabilité et les effets de ces accidents". De plus, un arrêté de 2005 (Journal Officiel, 2005) précise que "la probabilité peut être déterminée selon trois types de méthodes : de type qualitatif, semi-quantitatif ou quantitatif".

L'objectif de cette communication est de montrer certaines différences pour l'évaluation probabiliste des risques (EPR) technologiques selon l'approche : la première, semi-quantitative et plus rapide à mettre en œuvre, est basée sur des niveaux de confiance attribués aux barrières techniques de sécurité, tel que proposé dans le guide Oméga 10 (INERIS, 2008) ; la seconde, plus quantitative, est basée sur les arbres d'évènements et de défaillances.

## 2. Evaluation probabiliste des risques (EPR)

### 2.1. Introduction à l'EPR

L'évaluation probabiliste des risques (EPR) est une composante de l'analyse de risques, et plus spécifiquement de l'analyse quantitative des risques (AQR). L'EPR intervient après la phase d'identification des dangers et elle est ensuite un outil essentiel à la phase de gestion des risques (dont la réduction des risques fait partie). L'EPR est aujourd'hui un outil de maîtrise des risques largement utilisé et accepté dans de nombreux secteurs d'activité, dont le nucléaire, le spatial et l'aéronautique, l'industrie pétrolière et gazière, et l'industrie des procédés en général. L'EPR permet de caractériser chaque risque selon ses trois composantes, ce qui consiste à :

1. Identifier les scénarios possibles d'accident (i.e. les séquences d'évènements qui pourraient conduire à un accident, ou à tout autre évènement redouté) ;
2. Quantifier la fréquence de ces scénarios ;
3. Évaluer les conséquences de ces scénarios (i.e. les gravités).

Une EPR débute généralement par un évènement redouté central (ERC), pour lequel les causes sont décrites jusqu'aux événements initiateurs (EI) par un arbre de type "arbre de causes", et les conséquences sont décrites jusqu'aux évènements redoutés de fin (ERF) par un "arbre d'évènements". Ici, les ERF seront définis en termes de phénomènes dangereux (PhD), c'est-à-dire d'effets susceptibles d'entrainer des conséquences, et non par des accidents (i.e. qui tiennent aussi compte de l'exposition et de la vulnérabilité des cibles aux effets). Les "embranchements" des arbres (en amont et en aval de l'ERC) sont constitués des réalisations et non réalisations de fonctions de sécurité de prévention (FPrev., en amont de l'ERC) et de protection (FProt., en aval de l'ERC) réalisées par des barrières de sécurité, et d'évènements conditionnels (EC).

Chaque scénario (qui conduit à un PhD) est ensuite quantifié en termes de fréquence et de conséquences. L'évaluation des fréquences résulte d'une combinaison de fréquences (issues des EI) et de probabilités (issues des FPrev., FProt., et EC), qui doit respecter certaines règles mathématiques. Dans cette communication, deux approches sont étudiées pour l'évaluation des fréquences des PhD : semi-quantitative et quantitative. Celles-ci se différencient principalement par la méthode d'évaluation des probabilités de défaillance (i.e. non réalisations) des fonctions de sécurité.

### 2.2. Approche semi-quantitative pour l'évaluation des fonctions de sécurité

L'Institut National de l'Environnement Industriel et des Risques (INERIS) a développé une méthode semi-quantitative pour l'évaluation des fonctions de sécurité, basée sur les performances des barrières techniques de sécurité. Cette méthode, nommée "Oméga 10" (INERIS, 2008), est aujourd'hui couramment employée. Une méthode complémentaire, nommée "Oméga 20" (INERIS, 2009), concerne les barrières humaines de sécurité. Néanmoins, cette communication se focalisera essentiellement sur les barrières techniques de sécurité. Les performances de celles-ci sont évaluées selon trois critères :

- Efficacité (ou capacité de réalisation) ;
- Temps de réponse ;
- Niveau de confiance (NC).

Dans la suite, nous nous intéresserons plus spécifiquement au NC, qui est le critère quantitatif. Le NC est défini selon quatre niveaux allant de NC 1 à NC 4 et correspond à "la probabilité pour qu'une barrière technique de sécurité, dans son environnement d'utilisation, assure la fonction de sécurité pour laquelle elle a été choisie". Cette probabilité est calculée pour une capacité de réalisation et un temps de réponse donnés. De plus, le NC est lié à un facteur de réduction de risque (considéré comme l'inverse de la probabilité de défaillance de la fonction de sécurité) tel que défini dans la Table 1.

|      | Facteur de réduction de risque |
|------|-------------------------------|
| NC 1 | 10                            |
| NC 2 | 100                           |
| NC 3 | 1000                          |
| NC 4 | 10000                         |

**Table 1.** Niveau de confiance (NC) et facteur de réduction de risque

L'estimation du NC passe d'abord par l'évaluation qualitative de la barrière. Un certain nombre de critères est ainsi étudié, ce qui permet de juger de la performance de la barrière. Ensuite, un NC lui est attribué, en fonction de l'expérience, de l'état de l'art par type de dispositifs et des Tables 2 et 3. Ces tables sont une adaptation par l'INERIS de certaines exigences (liées aux contraintes architecturales) des normes de sécurité fonctionnelle CEI 61508 (CEI, 1998) et CEI 61511 (CEI, 2004). Ils fournissent la valeur maximale du NC en fonction de la tolérance aux anomalies du matériel (nombre de redondances) et de la proportion de défaillance en sécurité pour les éléments dits "simples" (i.e. sans microprocesseur) et "complexes" (i.e. avec microprocesseur). Les définitions suivantes sont applicables :

- la tolérance aux anomalies du matériel (HFT, pour "*hardware fault tolerance*") est de N si N+1 correspond au nombre minimal d'anomalies susceptibles de provoquer la perte de la fonction de sécurité ;
- la proportion de défaillance en sécurité (SFF, pour "*safe failure fraction*") est la proportion de défaillances (calculée à partir de taux de défaillance moyens) qui sont "en sécurité" ou "détectées en ligne" par rapport à l'ensemble des défaillances.

| Proportion de défaillance en sécurité (SFF) | Tolérance aux anomalies du matériel (HFT) | | |
|---|---|---|---|
| | 0 | 1 | 2 |
| < 60% | NC 1 | NC 2 | NC 3 |
| ≥ 60% à < 90% | NC 2 | NC 3 | NC 4 |
| ≥ 90% à < 99% | NC 3 | NC 4 | NC 4 |
| ≥ 99% | NC 3 | NC 4 | NC 4 |

**Table 2.** Niveau de confiance (NC) maximal pour des éléments "simples" (i.e. sans microprocesseurs)

| Proportion de défaillance en sécurité (SFF) | Tolérance aux anomalies du matériel (HFT) | | |
|---|---|---|---|
| | 0 | 1 | 2 |
| < 60% | . | NC 1 | NC 2 |
| ≥ 60% à < 90% | NC 1 | NC 2 | NC 3 |
| ≥ 90% à < 99% | NC 2 | NC 3 | NC 4 |
| ≥ 99% | NC 3 | NC 4 | NC 4 |

**Table 3.** Niveau de confiance (NC) maximal pour des éléments "complexes" (i.e. avec microprocesseurs)

### 2.3. Approche quantitative pour l'évaluation des fonctions de sécurité

L'approche quantitative pour l'évaluation des fonctions de sécurité est plus directement basée sur les normes de sécurité fonctionnelle CEI 61508 (CEI, 2010) et CEI 61511 (CEI, 2004). Ces normes concernent plus spécifiquement des systèmes instrumentés de sécurité (constitué de capteurs-transmetteurs, unités de traitement et actionneurs) mais sont généralisables à la plupart des autres barrières de sécurité. De plus, elles traitent de toutes les activités liées au cycle de vie des systèmes. Dans la suite, seule la "phase de réalisation" sera considérée. Celle-ci consiste à réaliser des systèmes en conformité avec la spécification des exigences de sécurité définis au préalable. De plus, en accord avec le sujet de cette communication, seules les exigences relatives à "la quantification de l'effet des défaillances aléatoires du matériel" seront considérées. En particulier, les exigences relatives aux "contraintes architecturales" et à l' "intégrité de sécurité systématique", par nature plus qualitatives que quantitatives, ne seront pas traitées ici.

La "quantification de l'effet des défaillances aléatoires du matériel" consiste à évaluer la probabilité moyenne d'une défaillance dangereuse lors de l'exécution sur sollicitation de la fonction de sécurité (PFDavg), dans le cas où la fonction de sécurité n'est réalisée que sur sollicitation et où la fréquence des sollicitations n'est pas supérieure à une par an (i.e. mode de fonctionnement à faible sollicitation). Ce résultat doit alors être inférieur ou égal à un "objectif chiffré de défaillance" fixé lors de la spécification des exigences de sécurité. Les exigences sur cette quantification consistent en une liste de paramètres à prendre en compte : l'architecture du système (en fonction de ses éléments) ; les taux de défaillance ; les causes communes de défaillance ; la couverture et les intervalles de temps des essais de diagnostic en ligne automatiques ; les intervalles de temps et l'efficacité des essais périodiques (essais destinés à détecter les défaillances "dangereuses" non détectées par les essais de diagnostic en ligne automatiques) ; les temps de réparation ; et l'effet d'erreurs humaines aléatoires. Certaines méthodes permettant cette quantification sont mentionnées, à titre informatif, dans la CEI 61508 (CEI, 2010) : équations simplifiées, diagrammes de fiabilité, arbres de défaillance, chaînes de Markov, réseaux de Petri. Une étude comparative a alors permis de conclure que ces différentes méthodes sont capables de produire des résultats similaires pour l'estimation de la PFDavg, en tenant compte des caractéristiques requises par la norme (Brissaud, 2012). Pour cela, il est néanmoins nécessaire de maîtriser la méthode utilisée et d'utiliser un outil logiciel adapté et performant. Dans cette communication, une approche par arbres de défaillance sera utilisée car son intégration dans un modèle d'EPR est plus intuitive.

## 3. Application de l'EPR selon les deux approches

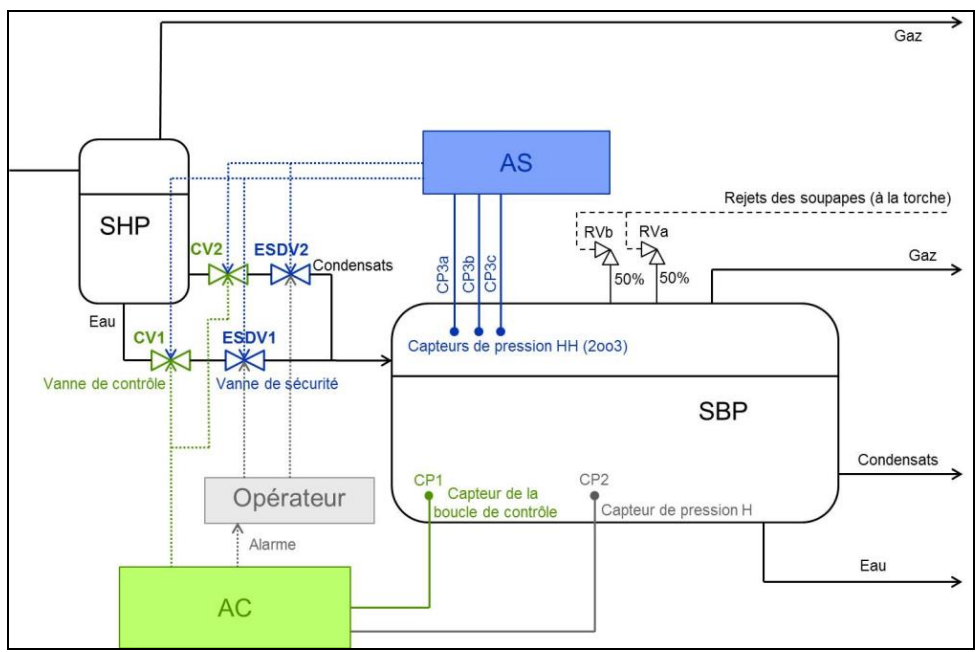

**Figure 1.** Schéma du cas d'étude

## 3.1. Cas d'étude

Le cas d'étude est un séparateur (gaz/condensats/eau) à basse-pression (SBP), installé en aval d'un séparateur à haute pression (SHP), tel que schématisé sur la Figure 1. Ces deux séparateurs sont reliés par deux lignes, une de condensats et une d'eau, dont le flux va du SHP vers le SBP. Sur chacune de ces lignes se trouvent une vanne de contrôle (CV) et une vanne de sécurité (ESDV). En sortie du SBP, du gaz, de l'eau et des condensats s'échappent par trois lignes distinctes. La pression dans le SBP est régulée via une boucle de contrôle composée d'un capteur de pression (CP1), d'un automate de conduite (AC) et des deux vannes de contrôle (CV1 et CV2).

L'évènement redouté central (ERC) considéré pour ce cas d'étude est un "dégagement de gaz suite à une surpression dans le séparateur basse-pression (SBP)". Afin de prévenir l'occurrence de cet ERC, les fonctions de sécurité de prévention sont réalisées par trois barrières techniques de sécurité :
1. Une alarme de pression haute avec action de l'opérateur, assurée par l'ensemble composé d'un capteur de pression (CP2), de l'automate de conduite (AC) et d'un opérateur agissant manuellement sur les deux vannes de sécurité (ESDV1 et ESDV2) ;
2. Un système instrumenté de sécurité (SIS), constitué de trois capteurs de pressions (CP3a, CP3b et CP3c) répondant à une logique en 2 sur 3 (le bon fonctionnement de 2 des 3 capteurs est suffisant pour assurer la fonction de sécurité), d'un automate de sécurité (AS) et des vannes de contrôle (CV1 et CV2) et de sécurité (ESDV1 et ESDV2) ;
3. Deux soupapes (RVa et RVb), sachant que le bon fonctionnement des deux soupapes est nécessaire pour assurer la fonction de sécurité (configuration en 2 fois 50%).

Selon les évènements initiateurs (EI) considérés, ces barrières de sécurité permettent ou non de prévenir l'occurrence de l'ERC, tel que décrit par l'arbre des causes de la Figure 2 (par exemple, le bon fonctionnement de n'importe laquelle de ces barrières est suffisant pour prévenir l'occurrence de l'ERC en présence d'EI 1, mais seul le bon fonctionnement des soupapes permet de prévenir l'occurrence de l'ERC en présence d'EI 4).

Si malgré tout l'ERC se produit, différents scénarios conduisant chacun à un phénomène dangereux (PhD) sont envisageables, selon l'occurrence ou non d'évènements conditionnels et de la réalisation ou non d'une fonction de sécurité de protection, tel que décrit par l'arbre d'évènements de la Figure 3.

Les EI sont listés dans la Table 4 avec leurs fréquences d'occurrence, données à titre d'exemple. Pour les barrières de sécurité, il sera considéré que seule une indisponibilité due à la non détection de défaillances "dangereuses" impliquera une non réalisation de la fonction de sécurité correspondante. Le taux de défaillance correspondant à ces défaillance est égale au taux de défaillance total de l'élément considéré, noté λ, que multiplie 1 moins la "proportion de défaillance en sécurité" (SFF). La période d'essais périodique permettant de détecter et de réparer ces défaillances est noté T1. Pour les vannes de sécurité (ESDV1 et ESDV2), des essais partiels sont également réalisés, selon une période T2, permettant de détecter et de réparer une proportion de ces défaillances notée PTC (pour "*partial test coverage*"). Enfin, un mode de défaillance commun à l'ensemble des capteurs (CP1, CP2, CP3a, CP3b et CP3c) est considéré, selon le modèle du facteur β (i.e. une proportion notée β des défaillances de chaque capteur implique une panne de l'ensemble des capteurs). Les valeurs de ces paramètres, adaptées de bases de données applicables au domaine pétrolier et gazier, sont reportées dans la Table 5. A noter que pour l'occurrence d'EI 1 ("Défaillance de la boucle de contrôle"), le taux de défaillance total des vannes de contrôle (CV1 et CV2) est considéré (i.e. utilisation de λ et non du produit (1 - SFF) × λ) car il s'agit ici d'un EI et non d'une fonction de sécurité. Enfin, l'action de l'opérateur en réponse à l'alarme de pression haute et la fonction de prévention réalisée par le local "blast proof" avec parois coupe-feu sont modélisés par les probabilités de défaillance à la sollicitation définies dans la Table 6, et les probabilités des évènements conditionnels sont définies dans la Table 7, données à titre d'exemple.

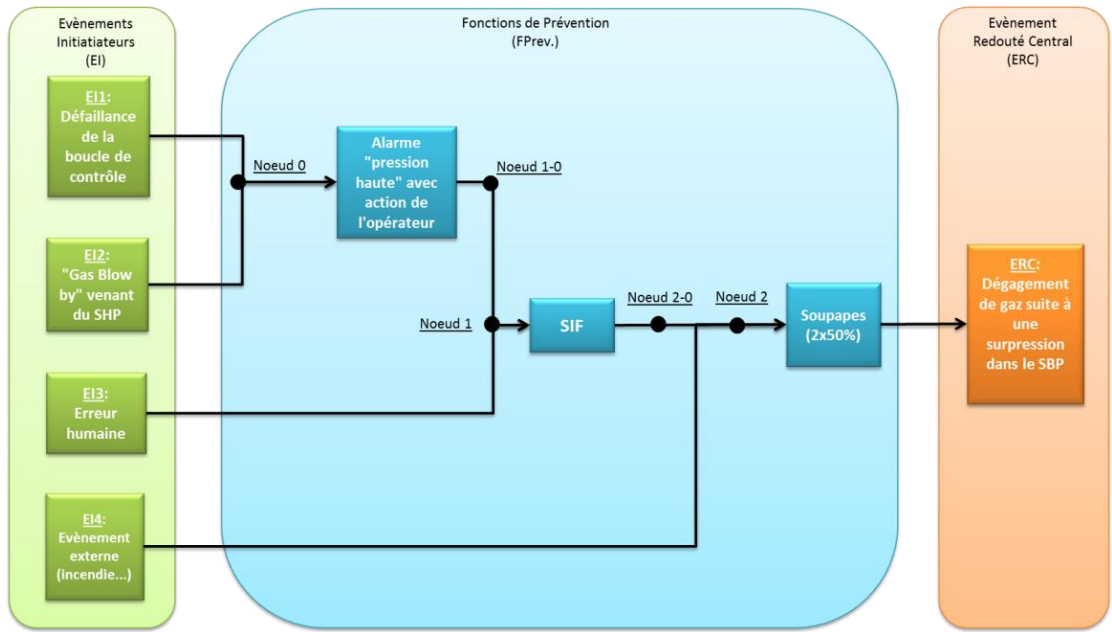

**Figure 2.** Modèle d'EPR : partie en amont de l'ERC

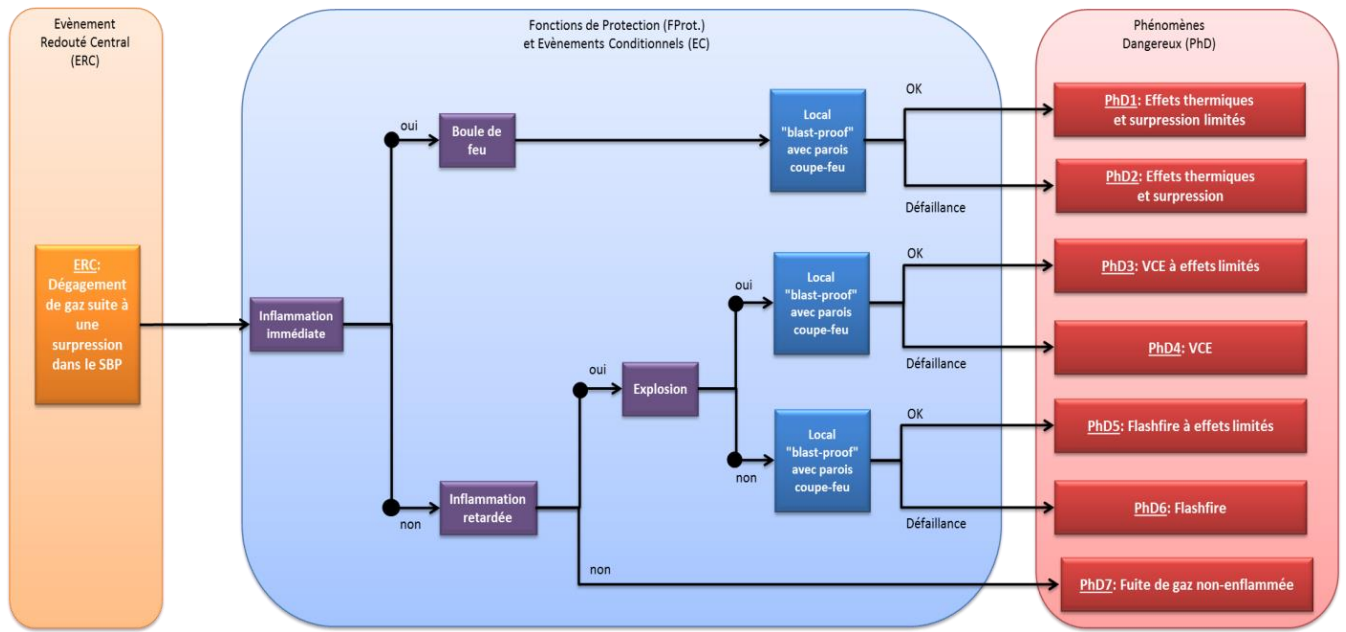

**Figure 3.** Modèle d'EPR : partie en aval de l'ERC

|  | Fréquence d'occurrence |
|---|---|
| EI 1 : Défaillance de la boucle de contrôle | Déterminée en fonction des paramètres de fiabilité et de maintenabilité des éléments de la boucle de contrôle |
| EI 2 : "Gas blow by" venant du SHP | 0,2 par an |
| EI 3 : Erreur humaine | 0,1 par an |
| EI 4 : Evènement externe (incendie…) | 0,005 par an |

**Table 4.** Evènements initiateurs (EI) avec fréquences d'occurrence

|  | $\lambda$ [par heure] | SFF | T1 [années] | Autres paramètres |
|---|---|---|---|---|
| Capteurs de pression (CP) | $3{,}200 \times 10^{-6}$ | 0,800 | 4 | $\beta = 0.05$ |
| Automate de conduite (AC) | $3{,}000 \times 10^{-6}$ | 0,835 | 4 | - |
| Automate de sécurité (AS) | $2{,}000 \times 10^{-6}$ | 0,950 | 4 | - |
| Vannes de contrôle (CV) | $3{,}260 \times 10^{-6}$ | 0,650 | 4 | - |
| Vannes de sécurité (ESDV) | $1{,}114 \times 10^{-5}$ | 0,625 | 4 | T2 = 6 mois, PTC = 0,900 |
| Soupapes (RV) | $1{,}392 \times 10^{-6}$ | 0,500 | 4 | - |

**Table 5.** Paramètres de fiabilité et de maintenabilité des éléments des barrières de sécurité

|  | Probabilité (constante) de défaillance à la sollicitation |
|---|---|
| Action de l'opérateur en réponse à l'alarme de pression haute | 0,1 |
| Fonction de prévention réalisée par le local "blast proof" avec parois coupe-feu | 0,01 |

**Table 6.** Probabilités de défaillance à la sollicitation

|  | Probabilité conditionnelle (constante) |
|---|---|
| Inflammation immédiate (sachant l'ERC) | 0,7 |
| Inflammation retardée (sachant l'ERC et la non inflammation immédiate) | 1,0 |
| Explosion (sachant l'inflammation retardée) | 0,4 |

**Table 7.** Probabilités des évènements conditionnels

### 3.2. Application de l'approche semi-quantitative

Les Tables 2 et 3 permettent de définir les niveaux de confiance (NC) des barrières de sécurité, en s'appuyant sur les SFF données dans la Table 5. Les paramètres définissant le NC sont basés sur les éléments les plus pénalisants. A noter que :
- L'alarme de pression haute avec action d'un opérateur est définie comme complexe à cause de l'automate de conduite (AC) ;
- Le SFF de l'alarme de pression haute avec action d'un opérateur ne prend pas en compte l'action de l'opérateur ;
- Le SIS est définie comme complexe à cause de l'automate de sécurité (AS) ;
- Le HFT du SIS est de 1 car malgré l'utilisation d'un seul automate de sécurité (AS), ce type d'automate inclut des redondances internes et, en ce qui concerne les capteurs de pression et les vannes du SIS, ceux-ci ont bien une HFT de 1 de par les redondances mise en place (CV et ESDV).

Ensuite, le facteur de réduction de risque est obtenu en fonction du NC, grâce à la Table 1. Les résultats d'évaluation des barrières de sécurité sont donnés dans la Table 8.

| Barrière de sécurité | Paramètres définissant le NC | | | Niveau de confiance (NC) | Facteur de réduction de risque |
|---|---|---|---|---|---|
| | complexité | HFT | SFF | | |
| Alarme de pression haute avec action d'un opérateur | complexe | 0 | ≥ 60% à < 90% | NC 1 | 10 |
| SIS | complexe | 1 | ≥ 60% à < 90% | NC 2 | 100 |
| Soupapes | simple | 0 | < 60% | NC 1 | 10 |

**Table 8.** Niveau de confiance (NC) des barrières de sécurité et facteurs de réduction de risque

Les fréquences d'occurrence des évènements initiateurs (EI) de la Table 4 (par défaut, la fréquence d'occurrence d'EI 1 est définie à 0,1 par an) sont ensuite propagées dans le modèle d'EPR par un système simple de divisions par les facteurs de réduction de risque, et d'additions.

### 3.3. Application de l'approche quantitative

Pour appliquer l'approche quantitative, le module Booléen de la suite GRIF (Satodev, 2014), développé par SATODEV et propriété de TOTAL, a été utilisé. La partie en aval de l'ERC est alors directement modélisée par un arbre d'évènements, tel que représenté sur la Figure 4. La partie en amont de l'ERC est quant à elle modélisée par un arbre de défaillances, dont un extrait est représenté sur la Figure 5. A noter que les éléments "EI1", "FPrev1", "FPrev2", et "FPrev3" de la Figure 4 sont eux-mêmes représentés par des arbres d'évènements, permettant notamment d'inclure les paramètres définis dans la Table 5.

Les PFDavg des barrières de sécurité sont obtenus grâce aux arbres de défaillance et sont reportés dans la Table 9. Cette table inclut également le facteur de réduction de risque, calculé par l'inverse de la valeur de PFDavg. Les fréquences d'occurrence de l'ERC et des phénomènes dangereux (PhD) sont ensuite obtenues directement par le logiciel, en appliquant les règles mathématiques de combinaison d'évènements et en associant notamment les fréquences d'occurrence des EI et les probabilités de défaillance des barrières de sécurité.

| Barrière de sécurité | PFDavg | Facteur moyen de réduction de risque |
|---|---|---|
| Alarme de pression haute avec action d'un opérateur | $1{,}45 \times 10^{-1}$ | 6,92 |
| SIS | $3{,}46 \times 10^{-3}$ | 288,86 |
| Soupapes | $2{,}40 \times 10^{-2}$ | 41,67 |

**Table 9.** PFDavg des barrières de sécurité et facteurs (moyens) de réduction de risque

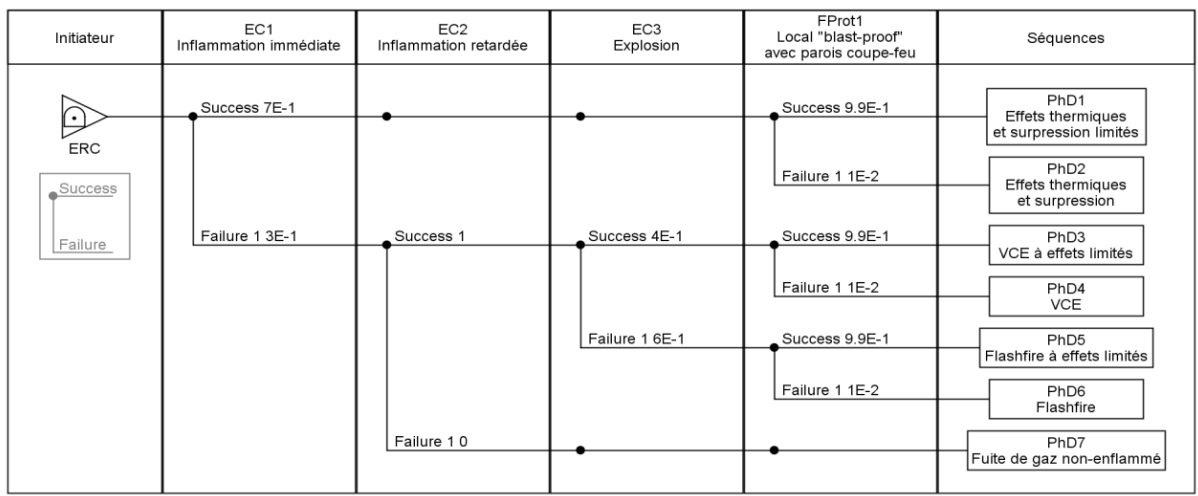

**Figure 4.** Arbre d'évènements obtenu avec le module Booléen de la suite GRIF

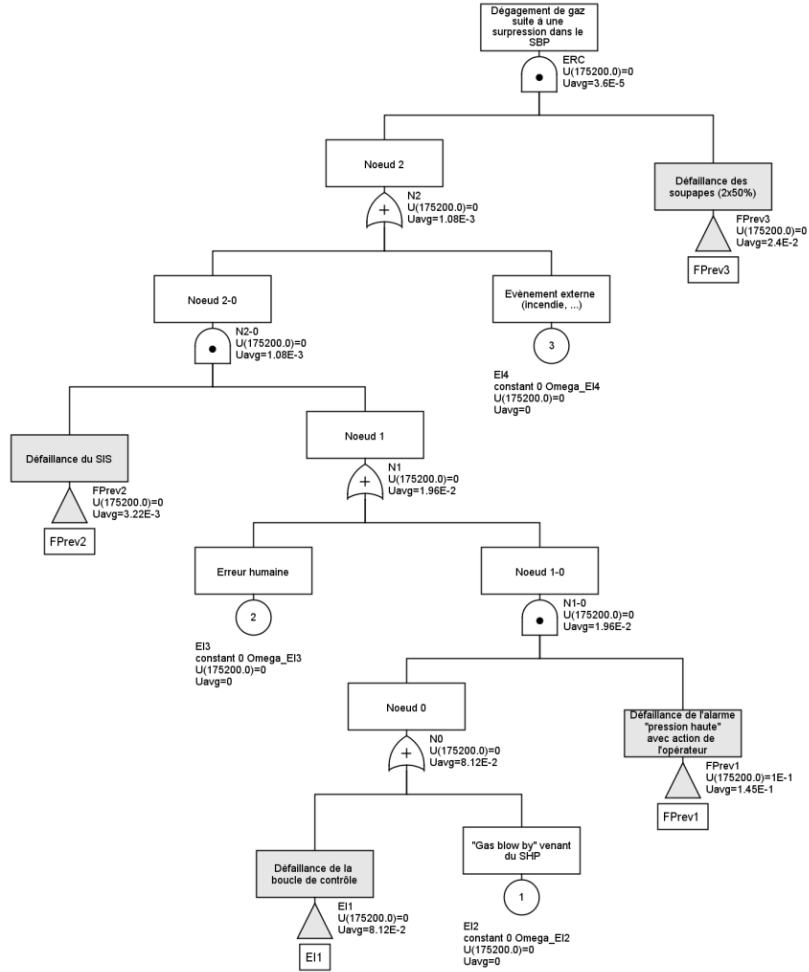

**Figure 5.** Extrait d'arbre de défaillances obtenu avec le module Booléen de la suite GRIF

## 4. Résultats et comparaisons des deux approches

### 4.1. Etudes de sensibilité

En plus du cas de base (Cas 0), une étude de sensibilité est proposée selon les quatre cas (Cas 1 à Cas 4) présentés dans la Table 10 ci-dessous.

| Cas | Description | Résumé |
|---|---|---|
| Cas 0 | Cas de base | Cas de base |
| Cas 1 | Multiplication de tous les taux de défaillance par 5 | $\lambda \times 5$ |
| Cas 2 | Division des proportions de défaillance en sécurité par 2 | SFF / 2 |
| Cas 3 | Multiplication par 2 des périodes d'essais (périodiques et partiels) | T1 x 2 et T2 x 2 |
| Cas 4 | Multiplication par 3 du facteur de mode commun de défaillance des capteurs | $\beta \times 3$ |

**Table 10.** Les différents cas de l'analyse de sensibilité

### 4.2. Résultats et comparaison des deux approches

Dans un premier temps, nous nous sommes intéressés aux fréquences annuelles d'occurrence de l'ERC ainsi que des PhD. Les résultats obtenus par les deux approches et les 5 cas considérés sont répartis dans la Table 11 ci-dessous.

On remarque notamment que, avec l'approche semi-quantitative, il n'y a que le Cas 2 (division des proportions de défaillance en sécurité par 2) qui présente des résultats différents. Les autres cas (Cas 1, 3 et 4) produisent les mêmes résultats que le cas de base (cas 0). Ceci est dû au fait que cette approche ne prend pas en compte certaines caractéristiques des barrières de sécurité telles que les taux de défaillance, la période des essais périodiques ou les modes communs de défaillance.

|       | Fréquences annuelles d'occurrence (moyenne sur 4 ans) | | | | | | | | | |
|-------|---|---|---|---|---|---|---|---|---|---|
|       | Approche quantitative | | | | | Approche semi-quantitative | | | | |
|       | **Cas 0** | Cas 1 | Cas 2 | Cas 3 | Cas 4 | **Cas 0** | Cas 1 | Cas 2 | Cas 3 | Cas 4 |
| ERC   | **2,07E-04** | 8,21E-03 | 7,79E-04 | 2,15E-04 | 2,33E-04 | **6,11E-04** | 6,11E-04 | 2,55E-03 | 6,11E-04 | 6,11E-04 |
| PhD1  | **1,43E-04** | 5,69E-03 | 5,40E-04 | 1,49E-04 | 1,61E-04 | **4,23E-04** | 4,23E-04 | 1,77E-03 | 4,23E-04 | 4,23E-04 |
| PhD2  | **1,45E-06** | 5,75E-05 | 5,46E-06 | 1,51E-06 | 1,63E-06 | **4,27E-06** | 4,27E-06 | 1,79E-05 | 4,27E-06 | 4,27E-06 |
| PhD3  | **2,46E-05** | 9,76E-04 | 9,26E-05 | 2,56E-05 | 2,77E-05 | **7,25E-05** | 7,25E-05 | 3,03E-04 | 7,25E-05 | 7,25E-05 |
| PhD4  | **2,48E-07** | 9,86E-06 | 9,35E-07 | 2,58E-07 | 2,79E-07 | **7,33E-07** | 7,33E-07 | 3,06E-06 | 7,33E-07 | 7,33E-07 |
| PhD5  | **3,69E-05** | 1,46E-03 | 1,39E-04 | 3,84E-05 | 4,15E-05 | **1,09E-04** | 1,09E-04 | 4,54E-04 | 1,09E-04 | 1,09E-04 |
| PhD6  | **3,73E-07** | 1,48E-05 | 1,40E-06 | 3,87E-07 | 4,19E-07 | **1,10E-06** | 1,10E-06 | 4,59E-06 | 1,10E-06 | 1,10E-06 |
| PhD7  | **0,00E+00** | 0,00E+00 | 0,00E+00 | 0,00E+00 | 0,00E+00 | **0,00E+00** | 0,00E+00 | 0,00E+00 | 0,00E+00 | 0,00E+00 |

**Table 11.** Fréquences annuelles de l'ERC et des PhD selon l'approche utilisée et le cas considéré

La Figure 6 présente les différences entre les résultats en termes de fréquence annuelle de l'ERC obtenus par l'approche quantitative et ceux obtenus par l'approche semi-quantitative. On observe ainsi certains écarts. Ceux-ci peuvent être dans un sens (par exemple, dans le Cas 1 l'approche quantitative fournit des résultats plus pessimistes d'un facteur 13 que l'approche semi-quantitative) ou dans un autre (par exemple, dans le Cas 2 l'approche quantitative fournit des résultats plus optimistes d'un facteur 3 que l'approche semi-quantitative).

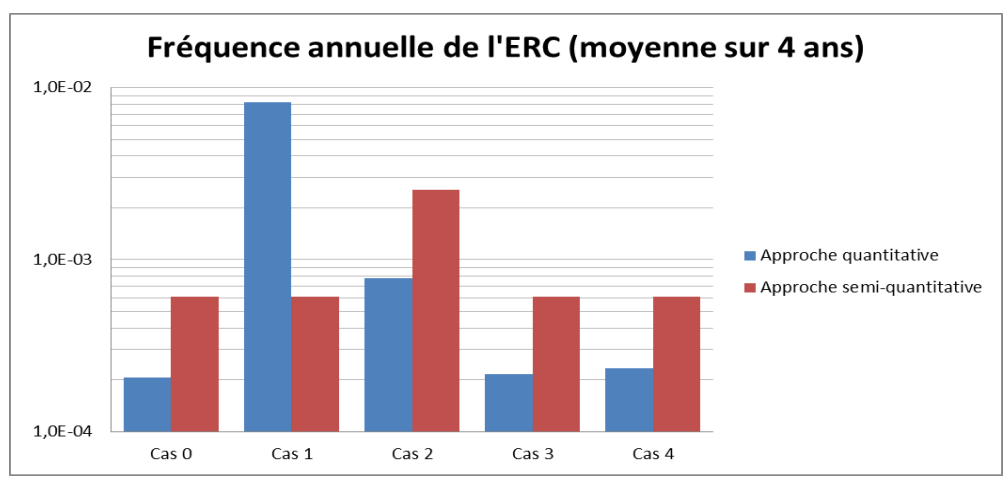

**Figure 6.** Comparaison de la fréquence d'occurrence de l'ERC selon l'approche utilisée

La Table 12 présente les facteurs de réduction du risque des trois barrières de sécurité de prévention selon l'approche utilisée et le cas considéré. On remarque que l'approche semi-quantitative ne fournit que des résultats de type "ordre de grandeur" (i.e. ordres de 10). Ces résultats ne sont cependant pas toujours plus pessimistes que ceux de l'approche quantitative. Par exemple, dans les cas 0, 1, 3 et 4, le facteur de réduction du risque de l'alarme de pression haute avec action de l'opérateur est ainsi plus élevé avec l'approche semi-quantitative qu'avec l'approche quantitative.

|          | Facteur de réduction du risque | | | | | | | | | |
|----------|---|---|---|---|---|---|---|---|---|---|
|          | Approche quantitative | | | | | Approche semi-quantitative | | | | |
|          | **Cas 0** | Cas 1 | Cas 2 | Cas 3 | Cas 4 | **Cas 0** | Cas 1 | Cas 2 | Cas 3 | Cas 4 |
| Alarme   | **6,92** | 3,36 | 4,96 | 6,32 | 6,92 | **10** | 10 | 1 | 10 | 10 |
| SIS      | **288,86** | 27,12 | 38,67 | 262,08 | 222,70 | **100** | 100 | 10 | 100 | 100 |
| Soupapes | **41,67** | 8,88 | 28,01 | 41,67 | 41,67 | **10** | 10 | 10 | 10 | 10 |

**Table 12.** Facteurs de réduction du risque des trois barrières de prévention selon l'approche utilisée et le cas considéré

## 5. Conclusion

En ce qui concerne la mise en œuvre de l'une et l'autre des approches, l'évaluation semi-quantitative selon le guide Oméga 10 est plus facile et plus rapide. En effet, le recours à une évaluation plus quantitative nécessite à la fois un outil logiciel adapté (l'utilisation de feuilles de calculs Excel ne permettant pas de prendre en compte convenablement les différentes dépendances au sein du modèle) et une certaine expertise spécifique.

L'évaluation des probabilités de défaillance des fonctions de sécurité (i.e. fonctions réalisées par les barrières de sécurité), selon l'approche semi-quantitative ne fournit que des résultats de type "ordre de grandeur" (i.e. ordres de 10). Pour autant, il est intéressant de constater que ces résultats ne sont pas toujours supérieurs ou inférieurs à ceux obtenus par l'approche quantitative. Pour le cas d'étude de base (Cas 0), la probabilité de défaillance d'une fonction de prévention (l'alarme) sur trois est ainsi plus faible avec l'approche semi-quantitative que l'approche quantitative.

Concernant l'évaluation des fréquences de l'ERC et des phénomènes dangereux du cas de base (cas 0), l'approche semi-quantitative a fourni des résultats qui sont plus élevés que ceux obtenus par l'approche probabiliste. Pour la fréquence de l'ERC, la différence est d'un facteur 3. Cependant, les études de sensibilité montrent aussi que, plusieurs caractéristiques (e.g. taux de défaillance et essais périodiques) n'étant pas être prises en compte dans l'approche semi-quantitative, les différences de résultats peuvent atteindre, jusqu'à un facteur 13 dans un sens (e.g. lorsque les valeurs des taux de défaillance ont été augmentées – cas 1), ou jusqu'à un facteur 3 dans un autre (e.g. lorsque les proportions de défaillance en sécurité ont été réduites – cas 2).

On en conclut ainsi que l'évaluation des risques peut être différente selon que l'approche utilisée soit semi-quantitative ou quantitative. Lorsque les enjeux sont importants, il convient donc de se questionner sur l'utilisation de l'une ou l'autre de ces approches proposées par la réglementation française sur la maîtrise des risques technologiques, afin d'opter pour la plus adaptées selon les informations disponibles et les besoins de l'étude.

## 6. Références

Brissaud F, Oliveira LF (2012), Average probability of a dangerous failure on demand: Different modelling methods, similar results, 11th International Probabilistic Safety Assessment and Management Conference and the Annual European Safety and Reliability Conference 2012 (PSAM11 ESREL 2012), Curran Associates, Inc., Red Hook (NY).

CEI (1998), CEI 61508 Ed. 1.0 : Sécurité fonctionnelle des systèmes électriques / électroniques / électroniques programmables relatifs à la sécurité, CEI, Genève.

CEI (2004), CEI 61511 : Sécurité fonctionnelle Systèmes instrumentés de sécurité pour le secteur des industries de transformation, CEI, Genève.

CEI (2010), CEI 61508-1 Ed. 2.0 : Sécurité fonctionnelle des systèmes électriques / électroniques / électroniques programmables relatifs à la sécurité, CEI, Genève.

INERIS (2008), Évaluation des performances des Barrières Techniques de Sécurité (DCE DRA-73) – Évaluation des Barrières Techniques de Sécurité - Ω 10, INERIS, Verneuil-en-Halatte.

INERIS (2009), Maîtrise des risques accidentels par les dispositions technologiques et organisationnelles (DRA 77) – Démarche d'évaluation des Barrières Humaines de Sécurité - Ω 20, INERIS, Verneuil-en-Halatte.

Journal Officiel (2000), Code de l'Environnement, Livre V : Prévention des pollutions, des risques et des nuisances (Partie législative), Titre I : Installations Classées pour la Protection de l'Environnement, JO n° 219 du 21 septembre 2000, Paris.

Journal Officiel (2003), Loi n°2003-699 du 30/07/03 relative à la prévention des risques technologiques et naturels et à la réparation des dommages, Chapitre II : Maîtrise de l'urbanisation autour des établissements industriels à risques, JO n° 175 du 31 juillet 2003, Paris.

Journal Officiel (2005), Arrêté du 29/09/05 relatif à l'évaluation et à la prise en compte de la probabilité d'occurrence, de la cinétique, de l'intensité des effets et de la gravité des conséquences des accidents potentiels dans les études de dangers des installations classées soumises à autorisation, JO n° 234 du 7 octobre 2005, Paris.

Lannoy A (2008), Maîtrise des risques et sûreté de fonctionnement – repères historiques et méthodologiques, Lavoisier, Paris.

Lenoble C, Durand C (2011), Introduction of frequency in France following the AZF accident, Journal of Loss Prevention in the Process Industries, 24 : 227-236.

Satodev (2014), http://www.satodev.com.